\documentclass[a4paper,11pt,DIV12]{scrartcl}

\usepackage[utf8]{inputenc}         %unter Linux
\usepackage{amsmath,amsfonts, amssymb, amsthm, mathrsfs}
\usepackage{bbm}
\usepackage{listings}
\usepackage{xcolor}
\usepackage{subfigure}
\usepackage{overpic}
\usepackage{url}

\usepackage{listings}  

\usepackage{tikz}
\usetikzlibrary{calc,trees,positioning,arrows,chains,shapes.geometric,%
    decorations.pathreplacing,decorations.pathmorphing,shapes,%
    matrix,shapes.symbols,backgrounds,shadows}

\title{A curved-element unstructured discontinuous Galerkin method on GPUs for the Euler equations}
\author{M. Siebenborn\thanks{Universität Trier, Universitätsring 15, D-54296 Trier, Germany, \texttt{(Siebenborn@uni-trier.de, Volker.Schulz@uni-trier.de)}} \and V. Schulz\footnotemark[1]  \and S. Schmidt\thanks{Imperial College London, Department of Aeronautics, South Kensington Campus, SW7 2AZ London, United Kingdom, \texttt{(s.schmidt@imperial.ac.uk)}}}

\date{}

\begin{document}
\maketitle

\begin{abstract}
In this work we consider Runge-Kutta discontinuous Galerkin methods (RKDG) for the solution of hyperbolic equations enabling high order discretization in space and time. We aim at an efficient implementation of DG for Euler equations on GPUs.
A mesh curvature approach is presented for the proper resolution of the domain boundary.
This approach is based on the linear elasticity equations and enables a boundary approximation with arbitrary, high order.
In order to demonstrate the performance of the boundary curvature a massively parallel solver on graphics processors is implemented and utilized for the solution of the Euler equations of gas-dynamics.
\end{abstract}

\begin{section}{Introduction}
The DG method is based on a discontinuous finite element spatial discretization originally introduced by Reed and Hill in the early 1970s for the neutron transport equation  \cite{reed_hill}.
Later, in a series of papers \cite{cockburn_shu_I,cockburn_shu_II,cockburn_shu_III,cockburn_shu_IV,cockburn_shu_V}, Cockburn and Shu combined the discontinuous Galerkin spatial discretization with an explicit Runge-Kutta time stepping (RKDG method) and extended the method to systems of conservation laws.
This created the opportunity for highly parallel implementations, since in the RKDG method, one grid cell only needs information from the immediate neighbouring cells to march in time.
Based on this, Biswas et al.\ investigated the potential of RKDG methods for parallelization in \cite{Biswas1994255}.\\
Another outstanding feature of DG methods is that the degree of basis functions can be chosen arbitrarily, thus leading to high order discretization.
However, attention has to be payed to the representation of curved boundaries.
Bassi and Rebay applied this method to the two dimensional Euler equations and worked out the importance of the boundary approximation with respect to the solution quality \cite{bassi_rebay_euler}.
They pointed out that the order of the method is limited by the order of the boundary representation.
Thus, it is necessary to deal with curved element discretizations, in order to overcome this issue.
While this seems not to be challenging for some test cases in two dimensions it is quite difficult to approximate complex geometries arising, e.g.\ , from industrial applications.\\
We present a mesh curvature procedure based on linear elasticity deformations which matches a desired boundary shape.
Similar techniques are well known for the deformation of discretization meshes and avoid expensive remeshing \cite{dwight2009robust}.
In order to deal with discontinuities arising in the solution of hyperbolic equations, slope limiters where introduced to the RKDG method by Cockburn and Shu.
Successful in the finite volume community, limiters are yet complicated for higher order methods.
This gives rise to the idea of adding artificial viscosity to the equations in order to smear out discontinuities and control the width of shocks.
Persson and Peraire proposed a detector to apply artificial viscosity only in the vicinity of shocks \cite{persson2006sub}.
This approach seems very promising because of its flexibility, since there is no dependency on the order of the scheme or the geometry of the discretization elements like for slope limiters.\\
Recently, there has been payed a lot of attention to discontinuous Galerkin methods because of their potential in terms of parallelization and high performance computing (HPC).
Especially the nodal DG method proposed by Hesthaven and Warburton \cite{hesthaven2002nodal} was shown by Klöckner, Warburton, Bridge and Hesthaven \cite{klockner2009nodal} to perform very efficiently on modern graphics processors (GPUs) gaining high speedups compared to conventional codes.\\
In this work, we present a novel approach combining the following aspects. We implement a high order Runge-Kutta discontinuous Galerkin method for the Euler equations of gas-dynamics on GPUs. For that, we propose an approach introducing unstructured, curved-element DG discretizations into a massively parallel GPU algorithm. Furthermore, we present a mesh curvature approach, which enables high order accurate boundary representations and seamlessly fits into the DG framework.
In particular, we cover the whole simulation chain from the generation of curved, body-fitted meshes up to the parallel HPC system solution. Finally, we demonstrate the performance of this algorithm on some challenging transsonic test cases, where discontinuities arise in the solutions.\\
This paper has the following structure. In section \ref{dg_method}, the discontinuous Galerkin method is shortly introduced and it is shown how the discrete operators are composed.
It is also covered how discontinuities in the solution are handled.
In section \ref{curvature}, we describe an approach for dealing with curvature in the DG discretization. Section \ref{gpu} shows the implementation on GPUs.
Finally, in section \ref{results} and \ref{conclusion} numerical results are presented and discussed.
\end{section}

\begin{section}{The discontinuous Galerkin method}\label{dg_method}
In this paper, we study a discontinuous Galerkin method for systems of hyperbolic conservation laws of the form
\begin{align*}
	\frac{\partial U}{\partial t} + \frac{\partial F_1(U)}{\partial x_1} + \dots + \frac{\partial F_d(U)}{\partial x_d} &= 0 & \text{in}\, (0,T] \times \Omega,\\
	U(0,x) &= U_0(x) & x \in \Omega
\end{align*}
where $U: \mathbb{R} \times \mathbb{R}^d \to \mathbb{R}^n \, , \quad U(t,\boldsymbol{x}) = \left(U_1(t,\boldsymbol{x}), \dots, U_n(t,\boldsymbol{x}) \right)^T$
is the vector of conserved quantities at a point $\boldsymbol{x}$ in $d$-dimensional space and at time $t$.
Here $\Omega \subset \mathbbm{R}^d$ is the domain of interest and $[0,T]$ a time interval. 
The vector fields $F_i: \mathbb{R}^n \to \mathbb{R}^n$ in this system are usually referred to as flux vectors.\\
For the purpose of this work we will maintain a more compact and widely used notation. Introducing the tensor
$F: \mathbb{R}^n \to \mathbb{R}^{n \times d} \, , \quad F(U) = \left( F_1(U)\; \dots \;  F_d(U) \right)$ in terms of the fluxes we then obtain
\begin{equation}\label{conservation_law}
	\frac{\partial U}{\partial t} + \nabla \cdot F(U) = 0.
\end{equation}
We follow the approaches in \cite{nodalDGbook} and mostly use the notation therein. 
For the derivation of the discontinuous Galerkin method, we assume that the domain of interest $\Omega$ is subdivided
into a finite set of $K$ disjoint, conforming elements $\Omega = \bigcup_{k=1}^{K} \Omega_k$. In our approach this is a tetrahedral mesh with curved elements.
The solution is then approximated using a space $\mathcal{V}_h$ of element-wise defined polynomials $\psi_j$ up to degree $p$
\begin{equation*}
\mathcal{V}_h = \bigoplus_{k=1}^K \mathcal{V}_h^k ,\quad  \mathcal{V}_h^k = \text{span} \{ \psi_j \left( \Omega_k \right), j=1,\dots, N_p \}.
\end{equation*}
Here, $N_p = \frac{(p+d)!}{d!p!}$ is the number of basis functions depending on the desired degree $p$ and the spatial dimension $d$.
Multiplying with a test function $\Phi$ from the same space and integration over the element $\Omega_k$ leads to
\begin{equation*}
\int\limits_{\Omega_k} \left[ \frac{\partial U}{\partial t} + \nabla \cdot F \left(U\right) \right] \Phi d \Omega = 0.
\end{equation*}
Integration by parts then yields the weak discontinuous Galerkin formulation
\begin{equation}\label{weak_dg}
\int\limits_{\Omega_k} \frac{\partial U}{\partial t} \Phi \, d \Omega = - \int\limits_{\partial\Omega_k} \left( F(U)^* \cdot \overrightarrow{n} \right) \Phi\, dS + \int\limits_{\Omega_k} F(U) \cdot \nabla \Phi d \Omega.
\end{equation}
Here, $\overrightarrow{n}$ denotes the outward pointing normal vector and $F^*$ an approximate Riemann solver, which deals with the double-valued state at the element interfaces, e.g.\ upwinding.
These approximate Riemann solvers are well known from finite volume context and a survey can be found in \cite{toro2009riemann}. 
For the purpose of this work, we tested both a local Lax-Friedrichs and a HLLC Riemann solver but without noticing major differences.
\begin{figure}
\begin{center}
	\begin{overpic}[width=1.0\textwidth]{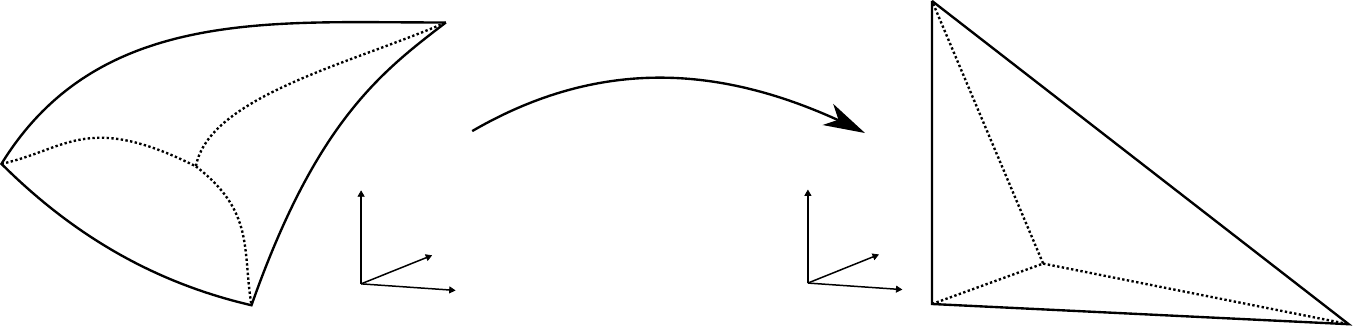}
	\put(44,20){$ \Psi(\boldsymbol x) = \boldsymbol r $}
	\put(28,10){$z$}
	\put(33.5,6){$y$}
	\put(34,0){$x$}
	\put(61,10){$t$}
	\put(66,6){$s$}
	\put(67,0){$r$}
	\end{overpic}
	\caption{Mapping from physical to reference element}
	\label{mapping}
	\end{center}
\end{figure}
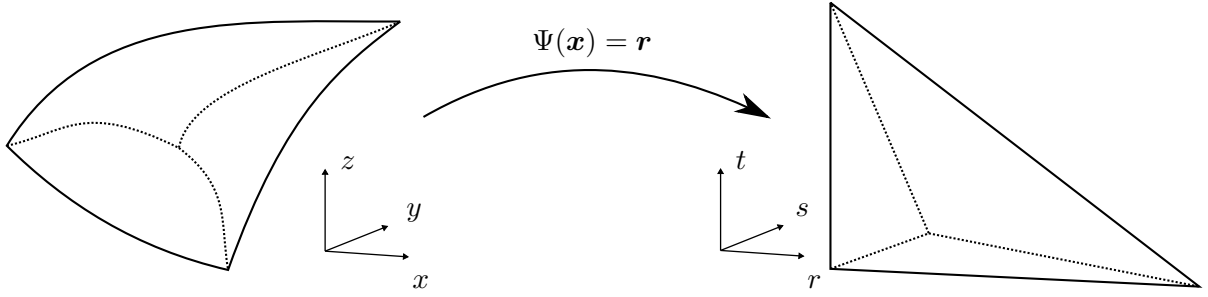\\
In the following, a discrete version of equation (\ref{weak_dg}) is derived, which can be implemented on a computer. In order to calculate the integrals occurring in (\ref{weak_dg}) we have to establish a link between an arbitrary curved element $\Omega_k$ in the mesh and the reference tetrahedron $\mathcal{T} = \lbrace -1 \leq r,s,t \leq 1\, , \, r+s+t \leq -1 \rbrace$ on which cubature/quadrature rules are known.
This is done by mapping functions $\Psi_k: \Omega_k \to \mathcal{T}$ for each element as illustrated in figure \ref{mapping} connecting the physical coordinates $\boldsymbol x = (x,y,z)$ to the computational ones $\boldsymbol r = (r,s,t)$.
Moreover, the partial derivatives of $\Psi_k$ have to be known since integration for any functions $f,g: \Omega_k \to \mathbb{R}^n$ in the physical space is evaluated as
\begin{equation*}
\int\limits_{\Omega_k} f g d\Omega = \int\limits_{T = \Psi_k\left(\Omega_k\right)} f\left(\Psi_k \right) g\left( \Psi_k \right) \left\vert \text{det} \left( D\Psi_k \right) \right\vert dT
\end{equation*}
with the convention
\begin{equation}
D\Psi = 
\begin{pmatrix}
r_x & r_y & r_z \\
s_x & s_y & s_z \\
t_x & t_y & t_z \\
\end{pmatrix}
\quad \text{and} \quad J = \left\vert \text{det} \left( D\Psi \right) \right\vert.
\end{equation}
Hence, it is sufficient to derive the local operators for the reference element and transform them into the physical space.
Then, we introduce a orthonormal, hierarchical set of modal basis functions $\lbrace\psi_i, i=1,\dots, N_p\rbrace$ on the reference element $\mathcal{T}$.
An arbitrary function $f$ on $\mathcal{T}$ is then interpolated on a given set of collocation points $\lbrace\boldsymbol r_i, i=1, \dots, N_p\rbrace$ (figure \ref{cubature_points}a) as
\begin{equation*}\label{modal}
f(\boldsymbol r_i) = \sum\limits_{j=0}^{N_p} \hat{f}_j \psi_j(\boldsymbol r_i) \, , \quad \forall i=1, \dots, N_p,
\end{equation*}
with modal expansion coefficients $\boldsymbol{\hat{f}} = (\hat{f}_1, \dots, \hat{f}_{N_p})^T$. The interpolation can be reformulated in terms of a multivariate Lagrange polynomial basis 
\begin{equation*}
f(\boldsymbol r_i) = \sum\limits_{j=0}^{N_p} f(\boldsymbol r_i) l_j(\boldsymbol r_i) \, , \quad \forall i=1, \dots, N_p
\end{equation*}
with nodal values $\boldsymbol{f} = (f(\boldsymbol r_1), \dots, f(\boldsymbol r_{N_p}))^T$.
These grid points are chosen according to \cite{warburton2006explicit} to ensure good interpolation properties.
Introducing the Vandermonde matrix $V_{ij} = \psi_j(\boldsymbol r_i)$ we obtain that these two formulations are linked as
\begin{equation*}
\boldsymbol{f} = V\boldsymbol{\hat{f}}.
\end{equation*}
Since there are as many basis polynomials as collocation points and depending on the interpolation quality $V$ is nonsingular.
\begin{figure}
\subfigure[Collocation points $\boldsymbol r_i$]{\includegraphics[width=0.33\textwidth]{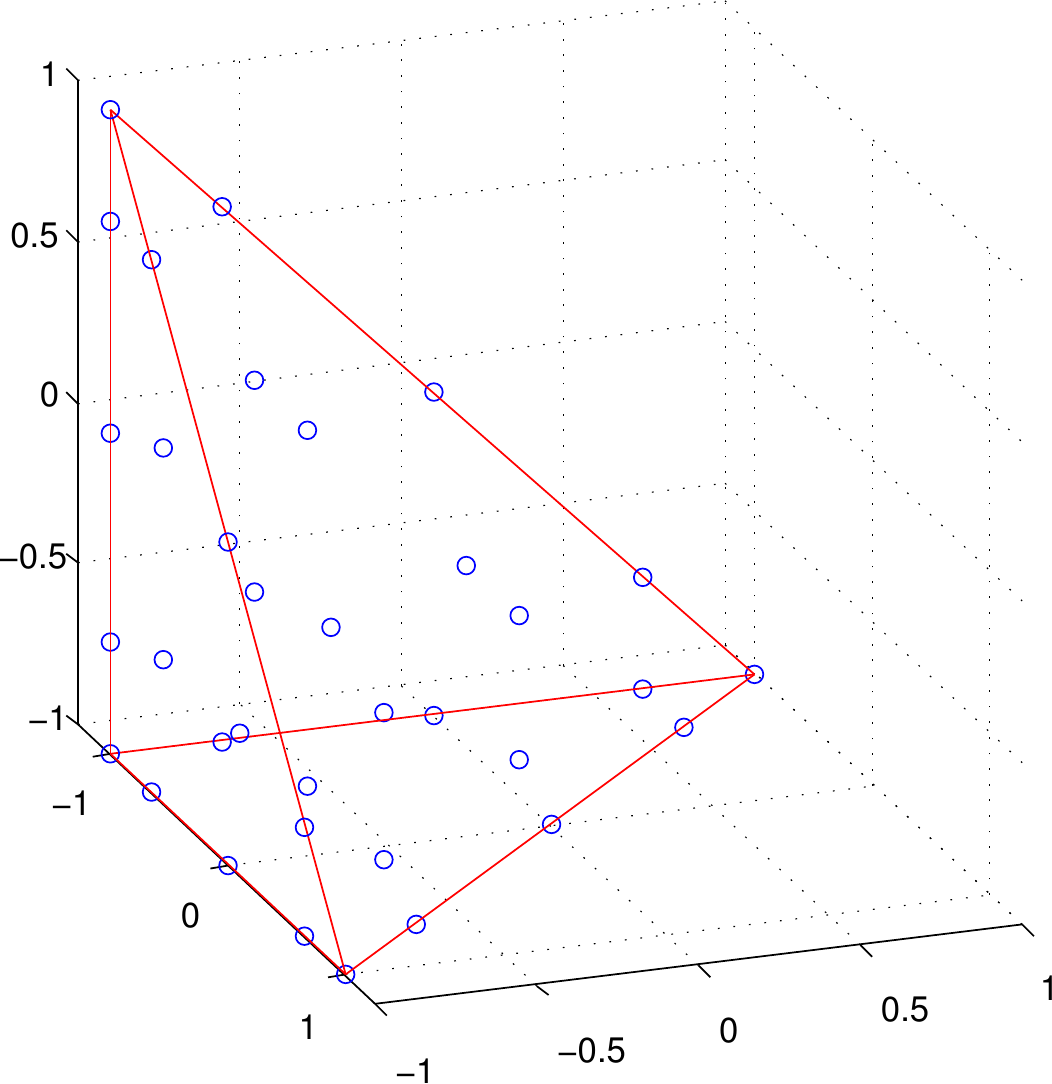}}\hfill
\subfigure[Volume cubature points  $\boldsymbol r_i^{\text{cub}}$]{\includegraphics[width=0.33\textwidth]{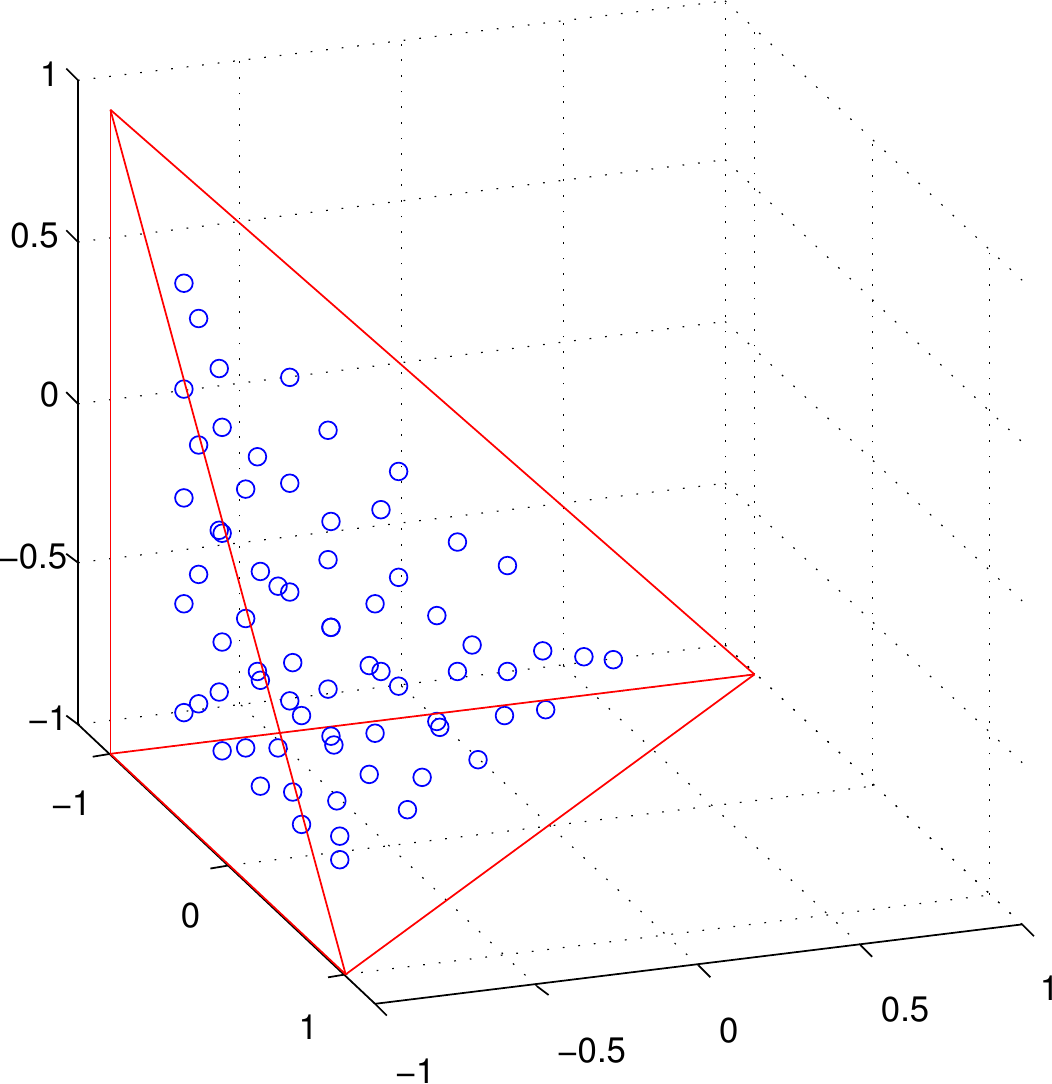}}\hfill
\subfigure[Surface quadrature points $\boldsymbol r_i^g$]{\includegraphics[width=0.33\textwidth]{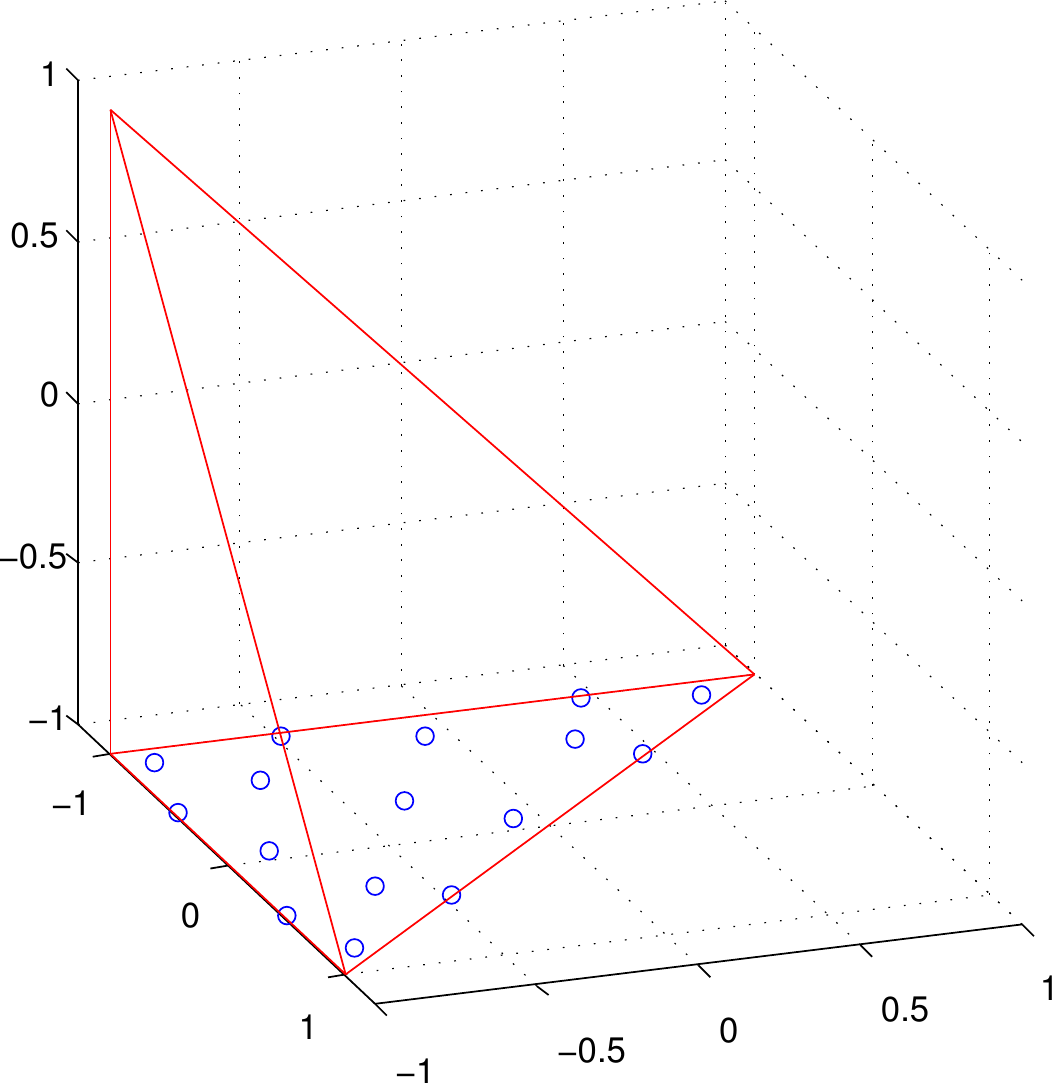}}
\caption{Collocation and quadrature points}
\label{cubature_points}
\end{figure}
Now, we are prepared to define the local operators for the integration in (\ref{weak_dg}).
For that purpose, we distinguish three sets of nodes: the collocation points $\boldsymbol r_i$ as introduced above and two sets of nodes used for integration.
These are on the one hand cubature points $\lbrace \boldsymbol r_i^{\text{cub}}, i=1, \dots,  N_{\text{cub}} \rbrace$ (figure \ref{cubature_points}b) for volume integrals \cite{grundmann1978invariant}
and on the other Gauss quadrature points $\lbrace \boldsymbol r_i^g, i=1, \dots,  N_g \rbrace$ (figure \ref{cubature_points}c) for surface integrals  \cite{cools1999cub}.\\
For the derivation of the discrete operators, we closely follow \cite{nodalDGbook}.
First, in order to interpolate the solution given at the nodal points $\boldsymbol r_i$ to quadrature points $\boldsymbol r_i^{\text{cub}}$ and $\boldsymbol r_i^g $ we need the following interpolation matrices
\begin{equation*}
\mathcal{I}_{\text{cub}} = V_{cub} V^{-1} \in \mathbb{R}^{N_{\text{cub}} \times N_p} \, , \quad \mathcal{I}_g = V_{g} V^{-1} \in \mathbb{R}^{N_{g} \times N_p}
\end{equation*}
where $\left(V_{\text{cub}}\right)_{ij} = \psi_j\left(\boldsymbol r_i^{\text{cub}}\right)$ and $\left(V_g\right)_{ij} = \psi_j(\boldsymbol r_i^g)$.
Multiplying the vector of nodal unknown values with these matrices first transforms them to modal expansion coefficients and then evaluates the modal basis functions at the cubature points.\\
For the construction of the stiffness matrices we need to evaluate derivatives like $\frac{\partial l_i}{\partial r}$ at the cubature points. Here we introduce the operator 
\begin{equation*}
D_r = V_r V^{-1} \in \mathbb{R}^{N_{\text{cub}} \times N_p}
\end{equation*}
where $(V_r)_{ij} = \left. \frac{\partial \psi_j(\boldsymbol r)}{\partial r} \right|_{\boldsymbol r = \boldsymbol r^{\text{cub}}_i}$.
The $s$ and $t$ derivatives are analogous.\\
With this we are now prepared to assemble the local stiffness matrix for one specific element leaving out the element number $k$ for clarity
\begin{equation}\label{stiffness}
S_x = \left( D_r^T \cdot \text{diag}(r_{x,i}) + D_s^T \cdot \text{diag}(s_{x,i}) + D_t^T \cdot \text{diag}(t_{x,i})\right) \cdot \text{diag}(J_iW_i^{\text{cub}}) \in \mathbb{R}^{N_p \times N_{\text{cub}}}
\end{equation}
where $W_i^{\text{cub}}$ are cubature weights. $S_y$ and $S_z$ are obtained in the same way.
The local mass matrices are then given by
\begin{equation}\label{mass}
M = \mathcal{I}_{\text{cub}}^T \cdot \text{diag}(J_iW_i^{\text{cub}}) \cdot \mathcal{I}_{\text{cub}} \in \mathbb{R}^{N_p \times N_p}
\end{equation}
Finally, we obtain the local face mass matrices
\begin{equation}\label{face_mass}
M_{\partial \Omega} = \mathcal{I}_{g}^T \cdot \text{diag}(J_iW_i^{g}) \in \mathbb{R}^{N_p \times N_g}
\end{equation}
where $W_i^{g}$ are the weights of a 2d Gauss integration rule for triangles.
Plugging these operators into equation (\ref{weak_dg}) yields the semi-discrete system 
\begin{equation}\label{semidiscrete}
\frac{\partial U_h}{\partial t} = M_k^{-1} \sum\limits_{m=1}^3 S_{k,x_m} F_m\left(U_h^{\text{cub}}\right) - M_k^{-1} M_{\partial \Omega_k} \left(F\left(U_h^{g}\right)^* \cdot \overrightarrow{n} \right).
\end{equation}
In the equation above, $U^{\text{cub}}_h = \mathcal{I}_{\text{cub}} U_h$ denote the degrees of freedom associated with the cubature points used for the volume integration (figure \ref{cubature_points}b)
and $U^{g}_h = \mathcal{I}_{g} U_h$ the ones for surface integration (figure \ref{cubature_points}c) respectively.
Since the operators are only local and therefore represented by small matrices, a LU or Cholesky decomposition of $M_k$ can be calculated in an initial step of the solver.\\
Finally, the system is discretized in time using a fourth-order accurate Runge-Kutta pseudo time stepping consisting of five stages.
We use a low storage version of this scheme as described in \cite{carpenter1994fourth} which does not require to keep the intermediate stages in memory.
Furthermore, the application of an explicit time integration like this RK method does not require to set up a global matrix system.
Thus, the method above is a so called matrix-free approach which is especially in three dimensions very attractive with respect to memory requirements.\\
When dealing with fluid dynamics (see section \ref{results}), the method described above works reliably for low Mach numbers.
Yet, for higher Mach numbers, shocks might appear in the solution leading to strong non-physical oscillations.
This phenomenon is well understood and can be overcome by adding a small amount of artificial viscosity to the equations as proposed by Persson and Peraire in \cite{persson2006sub}.\\
However, in regions where the solution is smooth there is no need for stabilization.
Hence, these authors compare the approximated solution of one component $u_h$ to the solution $\tilde{u}_h$ where they drop the modes with the highest frequencies
\begin{equation*}
u_h = \sum_{j=1}^{N_p} \hat{u}_j \psi_j \, , \quad \tilde{u}_h = \sum_{j=1}^{N_{p-1}} \hat{u}_j \psi_j
\end{equation*}
with the following smoothness indicator
\begin{equation*}
S_k = \frac{\int_{\Omega_k} \left( u_h - \tilde{u_h} \right) \cdot \left( u_h - \tilde{u_h} \right) d \Omega}{\int_{\Omega_k} u_h \cdot u_h \, d \Omega}.
\end{equation*}
In smooth regions of the solution, this indicator will be close to zero whereas it increases in regions with high frequencies.
The amount of artificial viscosity in element $k$ is then determined as
\begin{equation*}
\epsilon_k = \left\{
\begin{array}{l c l}
0 & \text{if} & s_k < s_0 -\kappa\\
\frac{\epsilon_0}{2}\left( 1 + \sin \frac{\pi \left( s_k - s_0\right)}{2\kappa}\right) & \text{if} & s_0 - \kappa \leq s_k \leq s_0 + \kappa\\
\epsilon_0 & \text{if} & s_k > s_0 + \kappa
\end{array}\right.
\end{equation*}
with $s_k = \log_{10} \left(S_k\right)$, $s_0 \sim \log_{10} \left(\frac{1}{p^4}\right)$ and empirically chosen parameters $\kappa$ and $\epsilon_0$.
Thus, by applying a shock detector  it can be decided where shocks arise and where to add viscosity.
The modified system of equation is then given by
\begin{equation*}
\frac{\partial U}{\partial t} + \nabla \cdot F(U) = \nabla \cdot \left( \epsilon \nabla U \right).
\end{equation*}
This can be rewritten as a system of first order equations and discretized using the same DG approach
\begin{align*}
\frac{\partial U}{\partial t} + \nabla \cdot F(U) -\nabla \cdot \sqrt{\epsilon} q &= 0\\
q - \sqrt{\epsilon} \nabla U & = 0.
\end{align*}
Through this procedure the semi-discrete system (\ref{semidiscrete}) is expanded by an additional equation and the viscous fluxes yielding (c.f.\  \cite{arnold2000discontinuous})
\begin{align*}
\frac{\partial U_h}{\partial t} &= M_k^{-1} \sum\limits_{m=1}^3 S_{k,x_m} \left[ F_m\left(U_h^{\text{cub}}\right) - \sqrt{\epsilon}q_h^{\text{cub}} \right] - M_k^{-1} M_{g} \left[ \left( F\left(U_h^{g}\right)^* - \left(\sqrt{\epsilon}q_h^{g}\right)^*\right) \cdot \overrightarrow{n} \right]\\
q_h &= \sum\limits_{m=1}^3 \sqrt{\epsilon} S_{k,x_m} U_h^{\text{cub}} - M_{g} \left[ \left(\sqrt{\epsilon}U_h^{g}\right)^* \cdot \overrightarrow{n} \right].
\end{align*}
\end{section}

\begin{section}{Mesh generation for DG}\label{curvature}
As pointed out in the introduction, a remarkable feature of the DG method is the freedom to choose the degree of basis functions.
However, because of that the mesh has to be dealt with caution.
It is tempting to reuse the meshes from the finite volume community. Yet, this leads to major problems, 
\begin{figure}[t]
\subfigure{\includegraphics[width=0.5\textwidth]{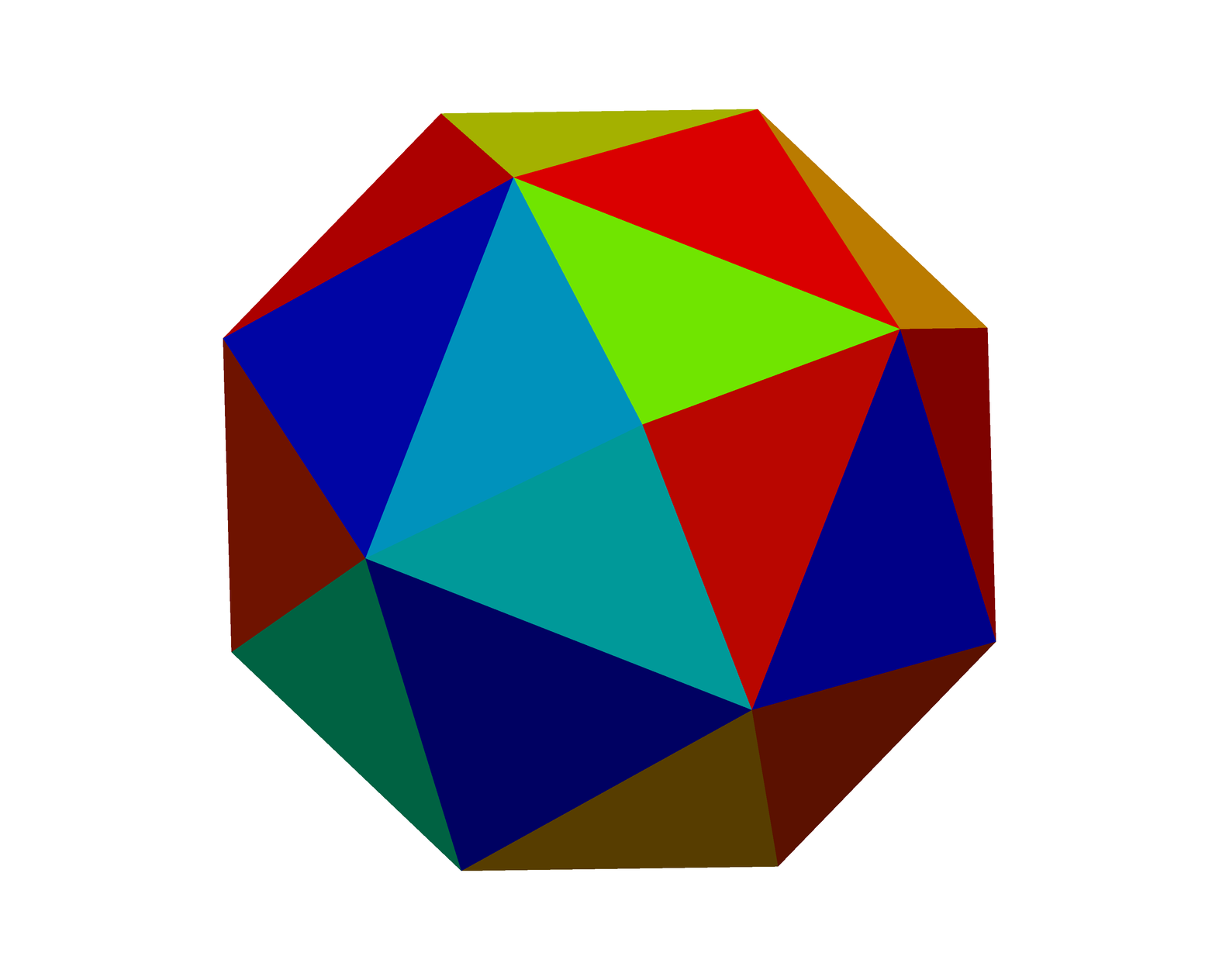}}\hfill
\subfigure{\includegraphics[width=0.5\textwidth]{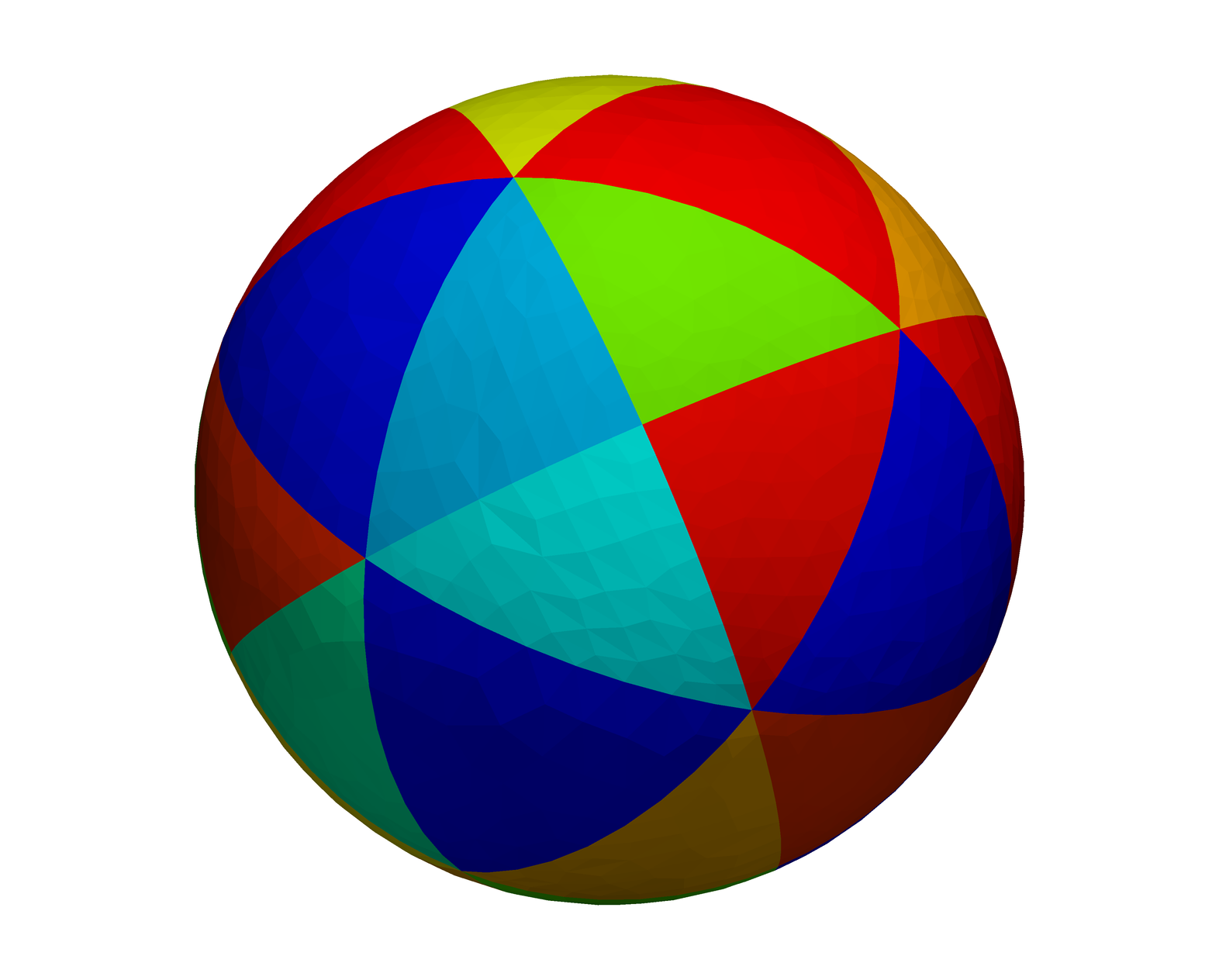}}
\caption{Deformation to a sphere}
\label{sphere}
\end{figure}
since the discretization with straight-sided elements leads to kinks at the boundary walls, which is problematic whenever the fluid is interacting with curved geometries.
In this situation the numerical solution might contain small non-physical shocks in each element at the boundary.
In principle, this could be overcome by using a very fine discretization. However, this is conceptually in contrast to higher order discretization methods.
Thus, it is necessary to introduce curved elements into the discontinuous Galerkin discretization in order to enable a higher order boundary representation.\\
An isoparametric mapping is applied for the realization of the element curvature. Assuming that the deformation of the elements can be represented in the same Lagrange polynomial basis as used for the DG scheme, we only need to know the displacement of the collocation points with respect to the ones in the straight sided element.
Thus, it is not enough to have a functional description of the curved boundary. Additionally, the locations of the collocation points have to be corrected such that they do not fall out of the element.
Moreover, this displacement procedure should work smoothly, in order not to perturb the integration.\\
In two dimensions this process is simple. Here, only the elements sharing a face with the boundary wall are affected.
By ensuring that the collocation points at the two other faces remain on their initial position, the neighboring elements remain untouched, even if they share one vertex with the boundary.
\begin{figure}[t]
\subfigure{\includegraphics[width=0.49\textwidth]{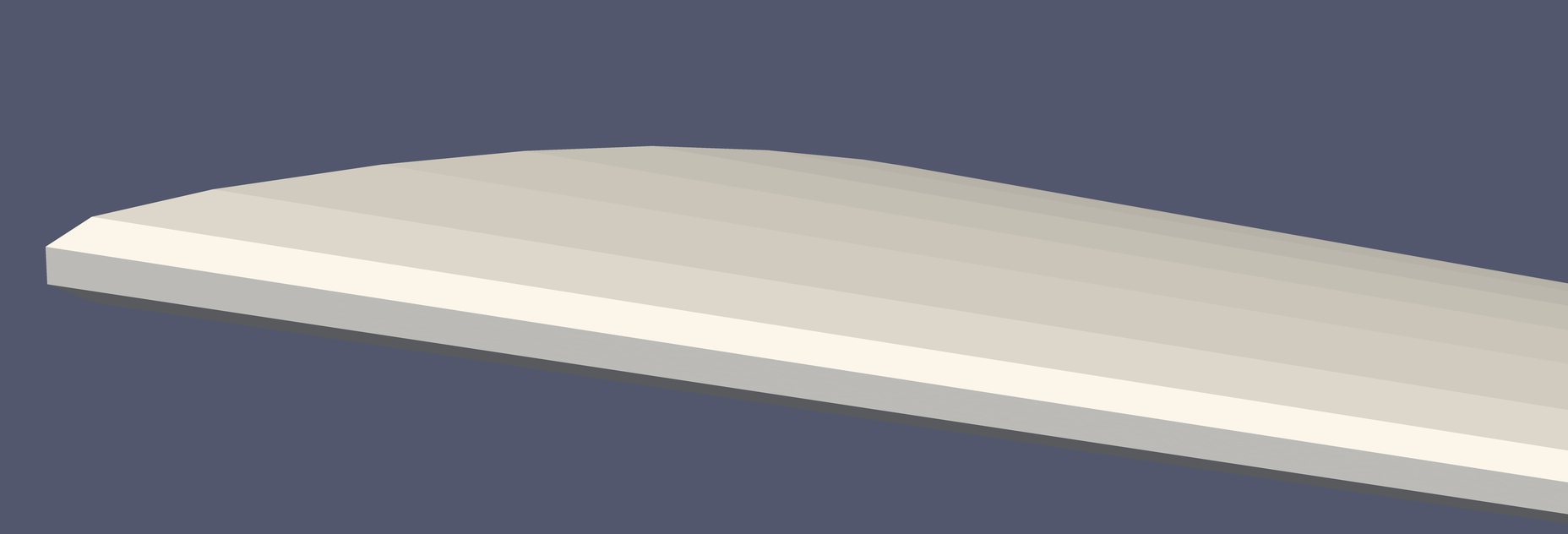}}\hfill
\subfigure{\includegraphics[width=0.49\textwidth]{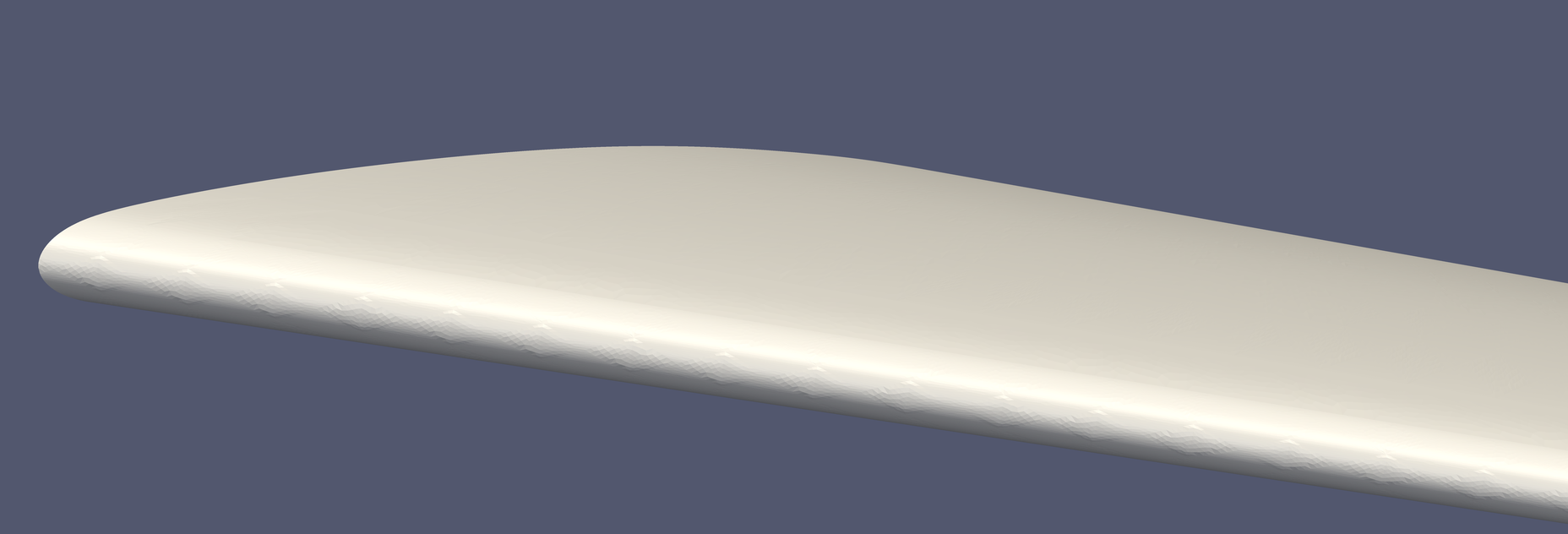}}
\caption{Leading edge of ONERA M6 airfoil surface - initial and curved geometry}
\label{onera}
\end{figure}
Obviously, this does not hold for the three dimensional case. By deforming a face on the boundary the other three faces in general have to be curved as well. This forces a whole layer of cells around the wall to be curved.\\
In the following, we focus on how to transport the curvature information of the boundary into the mesh. As mentioned before, in three dimensions a lot of cells are affected.
Hence, we suggest instead of tracking these cells to put a bounding box around the geometry marking an area in which cells get curved.
We follow the ideas in \cite{dwight2009robust} in order to model the discretization mesh as embedded in some flexible material and solve the linear elasticity equations with a finite element method.
By this, the curvature information of the boundary is transported into the mesh and can be retrieved at arbitrary points.\\
In order to keep this process as cheap as possible, we extract the cells which should become curved from the CFD-mesh and solve the linear elasticity equation on the resulting sub-mesh. 
The governing equations of the linear elasticity are given by
\begin{equation*}
\nabla \cdot \sigma = f \quad \text{on} \, \Omega_c.
\end{equation*}
Here $f$ denotes body forces acting on the solid, which will not be present in our case. $\Omega_c \subset \Omega$ is the region of the DG mesh containing the cells to be curved.
Finally, $\sigma$ denotes the stress tensor and is defined by the strain tensor $\epsilon$ as
\begin{equation*}
\sigma = \lambda \text{Tr}(\epsilon) I + 2\mu \epsilon \, , \quad \epsilon = \frac{1}{2} \left( \nabla u + \nabla u^T \right).
\end{equation*}
In the equations above, $\lambda$ and $\mu$ are the Lamé parameters, which can be expressed in terms of the Young's modulus $E$ and Poisson's ratio $\nu$
\begin{equation*}
\lambda = \frac{\nu E}{(1+\nu)(1-2\nu)} \, , \quad \mu = \frac{E}{2(1+\nu)} .
\end{equation*}
The vector-valued function $u: \Omega_c \to \mathbb{R}^3$ describes the displacement which the material would perform under the prescribed forces.
In our case, it is not important whether the chosen coefficients represent a physical material. Nevertheless, it is rather crucial to know how they act on the solution.
Young's modulus can be seen as a measure for the stiffness and Poisson's ratio describes how much a material expands in two coordinate directions when compressed in the third.
Since $f=0$ in our case the solution does not depend on the stiffness parameter $E$. Thus, the mesh deformation is controlled by $\nu$ only.
This parameter is then chosen depending on the situation. If we want to compute the flow past a sphere, like illustrated in figure \ref{curvature}, we apply a rubber-like material ($\nu \approx 0.5$) in order to obtain a smooth deformation, which propagates in each coordinate direction.
In the NACA0012 case, which is a 3d staggered version of the 2d test case, we choose $\nu$ close to $0$ in order to avoid deformation into the third dimension.\\
This system is then completed by the boundary conditions
\begin{align}\label{boundary}
u &= g \quad \text{on} \, \Gamma_{D1},\\
u &= 0 \quad \text{on} \, \Gamma_{D2},\\
\frac{\partial u}{\partial \overrightarrow{n}} &= 0 \quad \text{on} \, \Gamma_{N}.
\end{align}
We choose the boundaries such that $\Gamma_{D1}$ is the objects surface.
$\Gamma_{D2}$ is the bounding box defining the sub-mesh for the elasticity equation.
Finally, $\Gamma_N$ is used to model symmetry walls, where collocation points are allowed to slide but may not leave the plane.\\
One great advantage of this approach is the variable order of boundary representation.
Since, we are free to choose the degree of the finite element basis functions, we can fit the boundary representation to the order of the DG scheme.\\
Moreover, in most situations this approach avoids singularities due to the deformation.
This means, if the gap between the straight-sided mesh and the curved boundary is not to large, we can expect a feasible deformation of the mesh without overlappings or zero-volume cells.\\
However, this approach requires a parameterized representation of the surface, since the displacements at the boundary have to be evaluated at arbitrary coordinates.
While this seems not to be a problem for many geometries existing as CAD models, we are also dealing with geometries only available as a set of vertices.
Yet, there is a lot of literature dealing with the problem how to find a NURB surface representation to a given set of points.\\
For a parameterized surface it must be decided how the straight-sided mesh must be mapped to fit the curved geometry.
For this purpose a closest point problem has to be solved
\begin{equation}
\begin{aligned}
\min\limits_{(\alpha,\beta)} \quad& \frac{1}{2} \left\Vert \boldsymbol{x} - S(\alpha,\beta) \right\Vert^2\\
\text{s.t.} \quad& 0\leq \alpha,\beta \leq 1
\end{aligned}
\end{equation}
where  $S(\cdot,\cdot)$ denotes a NURB surface and $\boldsymbol x$ is a point on the straight-sided boundary. We solve this problem by the Gauss-Newton algorithm.
Let $(\alpha^*,\beta^*)$ be the optimal solution, then $g(\boldsymbol{x}) = S(\alpha^*,\beta^*) - \boldsymbol{x}$ is plugged into the boundary condition on $\Gamma_{D1}$.
\\With these boundary conditions we are now able to apply a high order finite element solver for the linear elasticity equation. For this purpose we use the GETFEM finite element toolbox.
Then, we store the solution obtained together with the CFD-mesh and use it later in the discontinuous Galerkin solver in order to retrieve mesh displacement information at arbitrary locations.

\begin{figure}[t]
\subfigure{\includegraphics[width=0.4\textwidth]{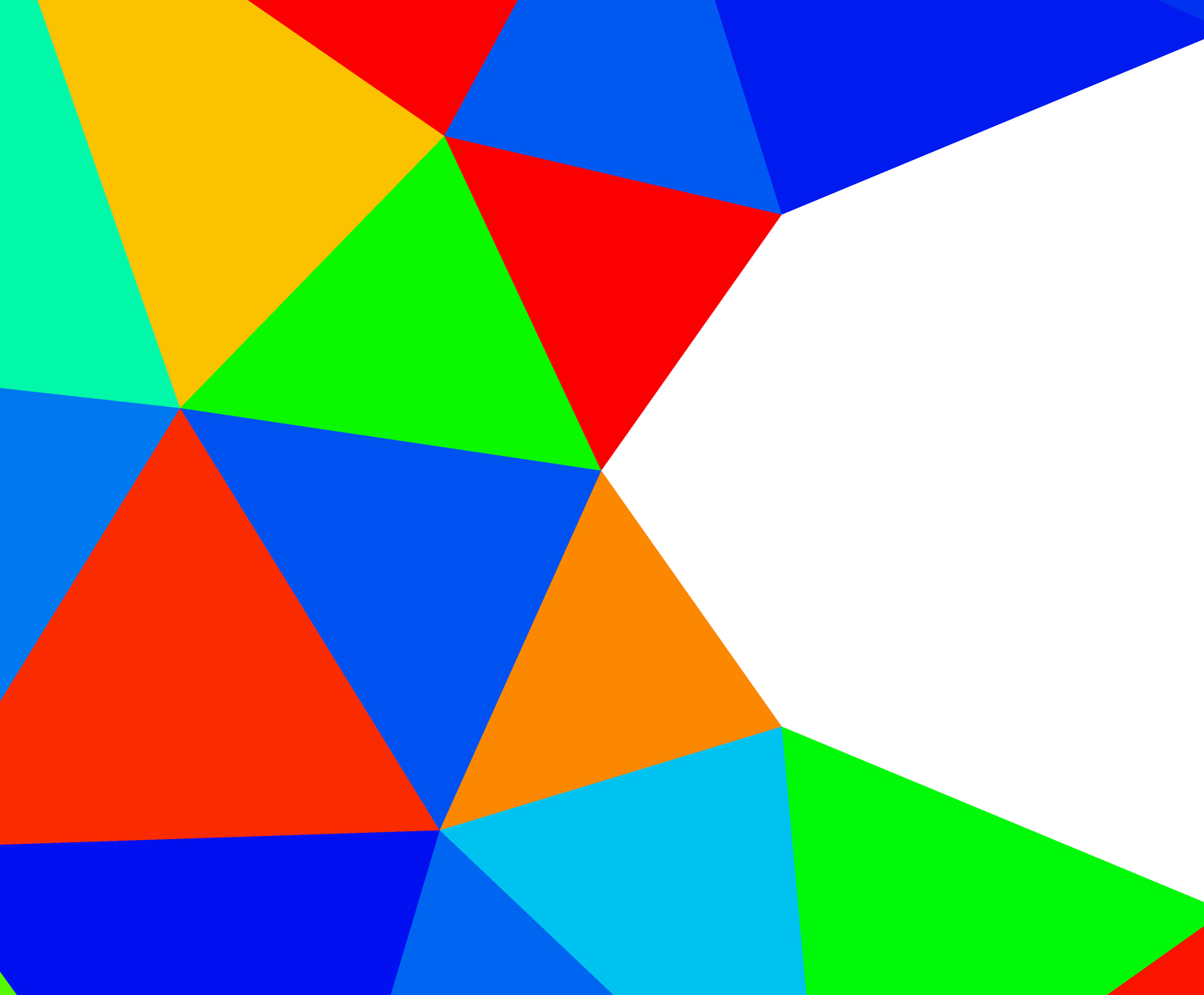}}\hfill
\subfigure{\includegraphics[width=0.4\textwidth]{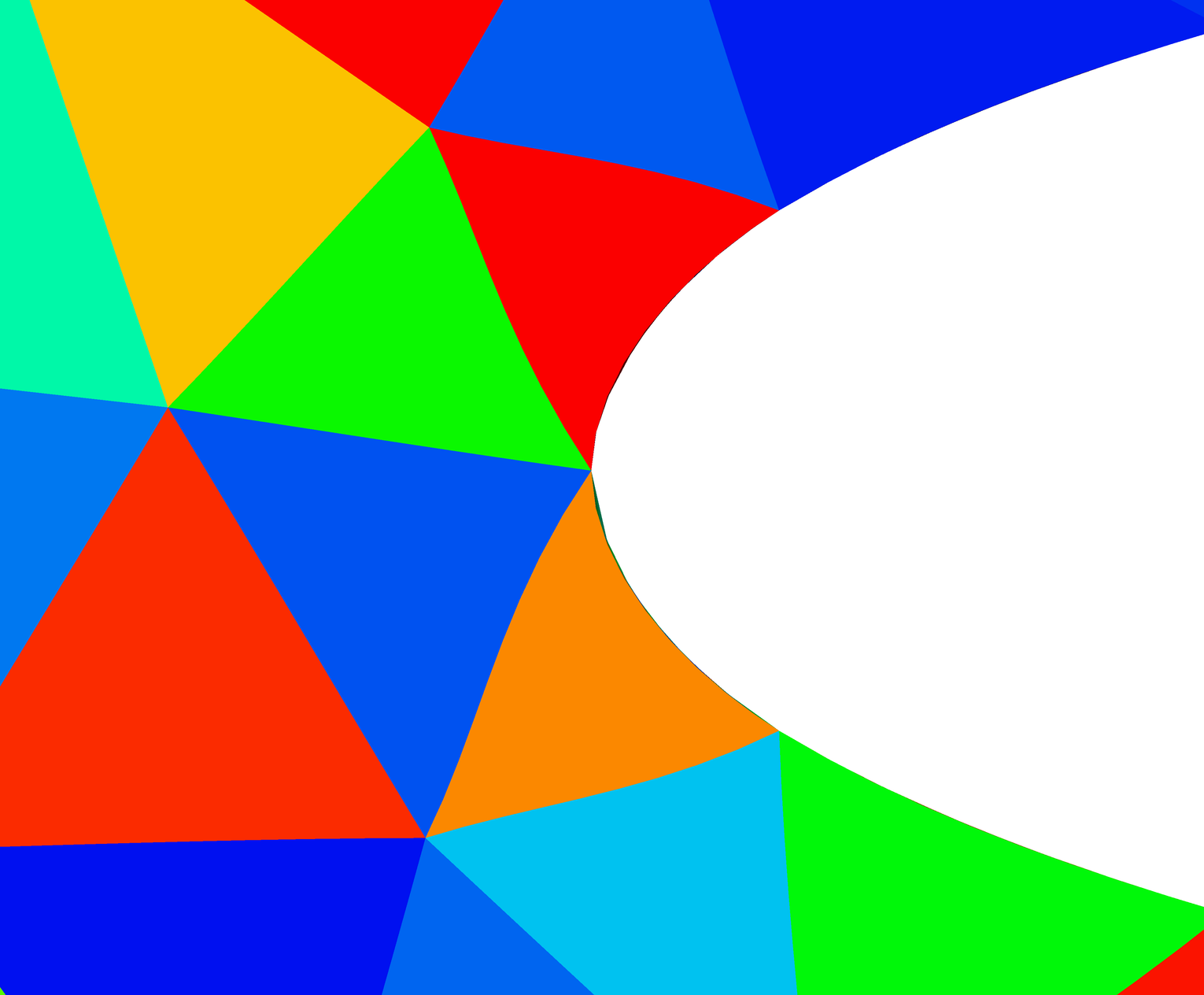}}
\caption{Tip of NACA0012 profile - initial and curved grid}
\label{naca}
\end{figure}

\end{section}

\begin{section}{Application on GPUs}\label{gpu}
Reviewing the DG method described above, it turns out that there is a lot of parallelism in the executions.
On the coarsest level one can proceed as in other discretization techniques like finite volume methods and partition the mesh.
Due to the nature of DG this does not lead to much data transfer between the compute nodes since only values at the element faces have to be communicated.\\
The same holds on a second level of parallelism. As seen in the derivation of the DG method in section \ref{dg_method}, there is only a loose coupling between elements.
Most operations can be performed independently. Just the evaluation of the surface integral requires the nodal values on the face of the neighbouring element to be known.\\
Furthermore, we can distinguish a third level of parallelism, the most interesting one.
Both on the volume and on the surface quadrature nodes, the nonlinear flux function $F$ has to be computed, which is a very expensive operation.
Not only the high order basis functions but also the rational Jacobian of the non-affine mapping force us to apply a very high order quadrature rule leading to a huge number of quadrature points. However, these operations are also independent per node.\\
These features suggest the application of graphics processing units, i.e.\ GPUs, to the DG method, since the design considerations above perfectly fit into the CUDA hardware model.
Originally designed to render many geometrical primitives or pixels in parallel the Nvidia Fermi architecture nowadays is based on a set of 16 so called streaming multiprocessors (SM).
Each of them features 32 cores which leads to a total number of 512 CUDA cores enabling hundreds of floating point operations at a time.\\
From the programming point of view, CUDA organizes threads in a hierarchical structure. On the lowest level there are the CUDA threads each having its own registers for temporary variables.
These threads are then gathered in thread blocks and identified through a three dimensional index. Each of these blocks comes with 64 kb shared memory in which data can be collected and distributed between the threads.
One major advantage of this concept is that shared memory is on chip.
Thus, data which is used multiple times or has to be shared between threads has to be requested from the global RAM of the GPU only once.\\
Finally, on the top level blocks are organized in a grid structure and again identified through three dimensional indices.\\
For further details we refer to the CUDA documentation \cite{nvidia2011}. Yet, one crucial aspect should be addressed.
In most cases the performance gain of a massively parallel GPU algorithm strongly depends on memory fetching strategies.
Generally, whenever the threads of a block $\lbrace t_0, \dots , t_n\rbrace$ fetch data from an array $A = \lbrace a_0, \dots, a_N \rbrace$ out of the global RAM of the GPU, these accesses should be organized in a blockwise fashion, where the blocksize is a multiple of $16$.
The threads $t_0, \dots , t_n$ should access successive data elements - each thread one element - where the first index is a multiple of $16$. Thus, memory accesses into the $i$th block should be organized as 
\begin{align*}
t_{\tau(j)} \leadsto a_{\left(16i + j \right)} \, , \quad j \in \lbrace 0,\dots,n\rbrace
\end{align*}
where $\tau$ is a permutation on $\lbrace 0,\dots,n\rbrace$.
In this case, multiple memory accesses can be performed at once leading to so called coalesced accesses.
Whereas, when a thread block requests data which is scattered throughout the memory this ends up in serialization of the accesses.\\
This issue is crucial particularly for older GPU hardware where finite volume methods (FVM) have been implemented for fluid dynamics.
Since the FVM can be seen as a special case of the DG method for $p=0$, there is only one unknown per cell. All neighbouring variables are involved in the update of the value of this unknown, which leads to a high ratio of scattered data accesses in case of unstructured meshes.
While GPU acceleration only seemed to work on structured meshes in FVM, the discontinuous Galerkin method overcomes this problem. 
Due to the higher order of the method, there are plenty of values per cell which can be organized successively in memory allowing high speedups even on unstructured meshes.\\
This hardware and programming layout can now be identified with the basic steps of the DG method.
As previously mentioned, operations on one element are mostly independent of the executions in other elements. Thus, it is obvious to assign one CUDA block to each DG element.
Then, the arithmetically expensive nonlinear flux evaluations can be performed in parallel.
Moreover, the subsequent matrix vector products can also be handled very efficiently within each block since we are dealing with small, dense matrices.
Finally, the CUDA execution handler decides how many thread blocks can be executed in parallel depending on the order of the DG method. \\
We have to reorganize the unknown values in order to meet the memory pattern of the GPU.
The length of each field should be a multiple of $16$ which is achieved by so-called zero padding as demonstrated in (\ref{padding}).
Here, $U_h^k$ denotes the unknown values inside element $k$.
As illustrated in chapter \ref{dg_method} we are dealing with a vector of $n$ unknowns interpolated at $N_p$ collocation points.
Thus, the padding length is determined such that the block length is enlarged from $N_p$ to $b_p = \lceil \frac{N_p +15}{16} \rceil$ and we ensure that the start address of each field is a multiple of $16$
\begin{equation}\label{padding}
\begin{matrix}
U_h^k =  &  \left[\vphantom{U^1_{N_p-1}}\right.&U^1_0 , & U^1_1 , & \dots , & U^1_{N_p-1} , & 0 , & \dots ,  & 0 ,&  \\
       &        &U^2_0 , & U^2_1 , & \dots , & U^2_{N_p-1} , & 0 , & \dots ,  & 0 ,&\\
       &        &\vdots  & \vdots  &         & \vdots        &\vdots &   & \vdots\vphantom{,}& \\
       &  &U^n_0 , & U^n_1 , & \dots , & U^2_{N_p-1} , & 0 , & \dots ,  & 0\vphantom{,} &\left. \vphantom{U^1_{N_p-1}} \right].\\
       &  &        &         &         &               & \multicolumn{3}{c}{\underbrace{\rule{1.8cm}{0pt}}_{\text{padding}}}  &
\end{matrix}
\end{equation}
The same works for the unknowns at the quadrature points.
In contrast, the values at the cubature points are only used once and do not need to get shared with other elements.
Thus, we do not need to store them in the global RAM.\\
In the following, we will consider the element specific volume integral in equation (\ref{weak_dg}) as an example.
For the integrations inside each element, we also need to fetch the element specific local operators from the global GPU RAM.
Thus, they are stored in a similar way like the vector of unknowns. For example the $S_x$ operator in equation (\ref{stiffness}) is stored as

\begin{equation}\label{padding2}
\begin{matrix}
S_x =   &  \left[\vphantom{U^1_{N_p-1}}\right.& s^x_{(0, 0)} , & s^x_{(1, 0)} , & \dots , & s^x_{(N_p-1,0)} , & 0 , & \dots ,  & 0 ,&  \\
       &  & s^x_{(0, 1)} , & s^x_{(1, 1)} , & \dots , & s^x_{(N_p-1,1)} , & 0 , & \dots ,  & 0 ,&\\
       &        &\vdots  & \vdots  &         & \vdots        &\vdots &   & \vdots\vphantom{,}& \\
       &  & s^x_{(0, N_{\text{cub}}-1)} , & s^x_{(1, N_{\text{cub}}-1)} , & \dots , & s^x_{(N_p-1,N_{\text{cub}}-1)} , & 0 , & \dots ,  & 0\vphantom{,} &\left. \vphantom{U^1_{N_p-1}} \right].\\
       &  &        &         &         &               & \multicolumn{3}{c}{\underbrace{\rule{1.8cm}{0pt}}_{\text{padding}}}  &
\end{matrix}
\end{equation}
In order to store the interpolation matrix $\mathcal{I}_{\text{cub}}$  the block length has to be as $b_{\text{cub}}  = \lceil \frac{N_{\text{cub}} +15}{16} \rceil $ leading to
\begin{equation}\label{padding2}
\begin{matrix}
\mathcal{I}_{\text{cub}} =   &  \left[\vphantom{U^1_{N_p-1}}\right.& \iota_{(0, 0)} , & \iota_{(1, 0)} , & \dots , & \iota_{(N_{\text{cub}}-1,0)} , & 0 , & \dots ,  & 0 ,&  \\
       &  & \iota_{(0, 1)} , & \iota_{(1, 1)} , & \dots , & \iota_{(N_{\text{cub}}-1,1)} , & 0 , & \dots ,  & 0 ,&\\
       &        &\vdots  & \vdots  &         & \vdots        &\vdots &   & \vdots\vphantom{,}& \\
       &  & \iota_{(0, N_p-1)} , & \iota_{(1, N_p-1)} , & \dots , & \iota_{(N_{\text{cub}}-1,N_p-1)} , & 0 , & \dots ,  & 0\vphantom{,} &\left. \vphantom{U^1_{N_p-1}} \right].\\
       &  &        &         &         &               & \multicolumn{3}{c}{\underbrace{\rule{1.8cm}{0pt}}_{\text{padding}}}  &
\end{matrix}
\end{equation}
\tikzset{mybox/.style={
rectangle,
rounded corners=2mm,
thick,
draw=black,
text width=15em,
text centered,
drop shadow,fill=white}
}
\begin{figure}[t]
\begin{center}
\begin{tikzpicture}[node distance=.5cm,>=stealth',bend angle=45,auto]
\tikzstyle{every node}=[font=\small]
\node[mybox] (mesh) {Mesh partitioning};
\node[mybox] at ($(mesh.south) - (0,2.5)$) (curve) {Curve mesh according to precomputed linear elasticity solution};
\path (mesh) edge[->,thick] (curve);
\node[mybox, below = of curve] (setup) {Setup operators and transfer to GPU};
\path (curve) edge[->,thick] (setup);
\node[mybox] at ($(setup.south) - (0,1.7)$) (interp) {Interpolate values to quadrature nodes};
\path (setup) edge[->,thick] (interp);
\node[mybox, below = of interp] (download) {Distribute nodes on processor cuts to adjacent processes};
\path (interp) edge[->,thick] (download);
%\draw [<->,thick] ($ (download.east) $) -- ($ (download.east) + (2cm,0) $) ;
%\draw [<->,thick] ($ (download.west) $) -- ($ (download.west) - (2cm,0) $) ;
%\draw [<->,thick] ($ (download.east) + (4cm,0)$) -- ($ (download.east) + (6cm,0) $) ;
%\draw [<->,thick] ($ (download.west) $) -- ($ (download.west) - (2cm,0) $) ;
\node[mybox, below = of download] (vol) {Interpolation to cubature nodes and volume integration};
\path (download) edge[->,thick] (vol);
\node[mybox, below = of vol] (surf) {Surface integration};
\path (vol) edge[->,thick] (surf);
\node[mybox, below = of surf] (RK) {Evaluate Runge-Kutta stage};
\path (surf) edge[->,thick] (RK);
\draw [->,thick] ($ (RK.west) $) -- ++(-0.8,0) -| ($ (interp.west) - (0.8,0) $) -- ($ (interp.west) $) ;

\node[mybox,text width=7em] at ($(download.west) - (3.0,0)$) (other1) {Communication with other compute nodes};
\node[mybox,text width=7em] at ($(download.east) + (2.5,0)$) (other2) {Communication with other compute nodes};
\path (download) edge[<->,thick] (other1);
\path (download) edge[<->,thick] (other2);

%\node[mybox] at ($(download.west) - (2,0)$) (receive1) {Receive face values};

\node[rectangle,rounded corners=2mm,thick,draw=black,fill=red!50] at ($(curve.west) - (0,-1.5)$) {Node $j$};
\node[rectangle,rounded corners=2mm,thick,draw=black,fill=red!50] at ($(interp.west) - (-0.3,-1.0)$) {GPU $j$};

%\node[rectangle,rounded corners=2mm,thick,draw=black,fill=red!50] at ($(curve.west) - (3.7,-1.5)$) {Node $j-1$};
%\node[rectangle,rounded corners=2mm,thick,draw=black,fill=red!50] at ($(interp.west) - (3.4,-1.5)$) {GPU $j-1$};

%\node[rectangle,rounded corners=2mm,thick,draw=black,fill=red!50] at ($(curve.west) + (8.8,+1.5)$) {Node $N$};
%\node[rectangle,rounded corners=2mm,thick,draw=black,fill=red!50] at ($(interp.west) + (9.1,+1.5)$) {GPU $N$};

\begin{pgfonlayer}{background}
	\filldraw [rounded corners=2mm,black!10] ($(RK.south) +(-4.0,-0.6)$) rectangle ($(curve.north) +(3.5,1.4)$);
	\filldraw [rounded corners=2mm,yellow!30] ($(RK.south) +(-3.6,-0.4)$) rectangle ($(interp.north) +(3.2,1.0)$);

%	\filldraw [rounded corners=2mm,black!10] ($(RK.south) +(4.5,-0.6)$) rectangle ($(curve.north) +(7.5,1.4)$);
%	\filldraw [rounded corners=2mm,yellow!30] ($(RK.south) +(4.9,-0.4)$) rectangle ($(interp.north) +(7.3,1.4)$);

%	\filldraw [rounded corners=2mm,black!10] ($(RK.south) +(-8.0,-0.6)$) rectangle ($(curve.north) +(-5,1.4)$);
%	\filldraw [rounded corners=2mm,yellow!30] ($(RK.south) +(-7.6,-0.4)$) rectangle ($(interp.north) +(-5.3,1.4)$);
\end{pgfonlayer}

\end{tikzpicture}
\caption{Algorithm execution on one compute node}
\label{flowchart}
\end{center}
\end{figure}
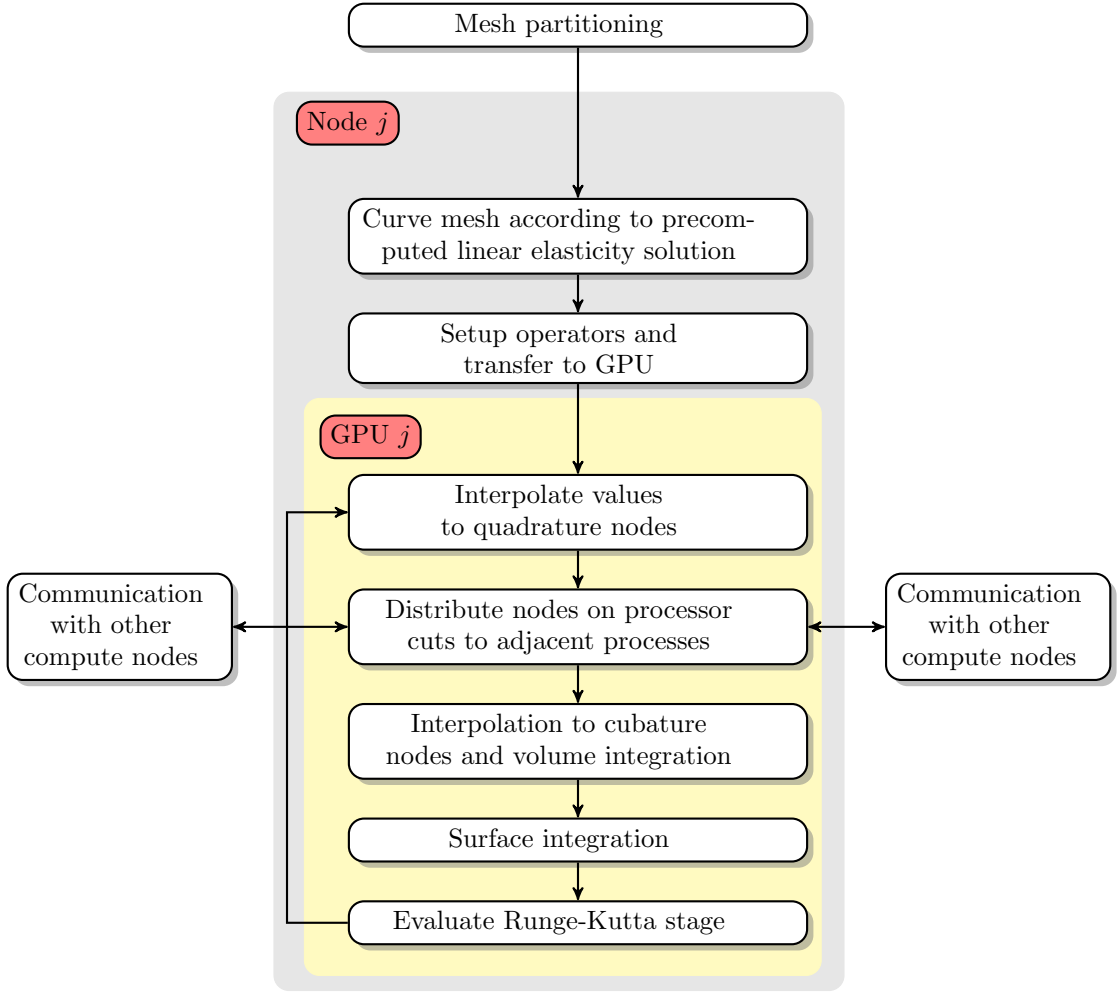\\
Note that these operators are stored in column-major form, which is important for the matrix-vector product.\\
With this data structures in global GPU RAM we implement the CUDA grid to have $K$ thread blocks, which will be associated with the DG elements.
Each of these blocks should contain $N_{\text{cub}}$ threads $\lbrace t_0, \dots, t_{N_{\text{cub}}-1}\rbrace$ and a matrix $U[n][N_\text{cub}]$ in shared memory to cache the unknown values.
At this point, we can assume that the matrix $U$ is large enough to store both $U_h$ and $U_{\text{cub}}$ since $N_{\text{cub}} > N_p$.
This holds, because for a proper integration of the nonlinear flux function there have to be more cubature points than collocation points.
The initial step is to cache the nodal values $U_h$ into the shared matrix. For that purpose the first $N_p$ threads load and store the $n$ fields of $U_h$ into $U$.
Then the product $\mathcal{I}_{\text{cub}} U^k_h$ is evaluated column-wise. The threads $t_0, \dots, t_{N_{\text{cub}}-1}$ load one column of $\mathcal{I}_{\text{cub}}$, each thread one value.
Since $\mathcal{I}_{\text{cub}}$ is stored in column-major form, these loads are coalescent.
Then the $n$ values in $U$ corresponding to this column are broadcasted to all threads and the multiplication is performed.
The broadcast operation is very efficient, since the matrix $U$ was fetched to the shared memory which is on chip in contrast to the global GPU RAM.
This is done successively for all columns and each thread is tracking the sum over all columns in its registers.
Finally, by this procedure, we obtain the nodal values interpolated to the cubature points $U_h^{\text{cub}}$, which are again stored in $U$.
Through this procedure the dot products between the $N_{\text{cub}}$ rows of $\mathcal{I}_{\text{cub}}$ and the fields of $U$ are performed in parallel each by one CUDA thread.\\
In order to calculate the volume integral in (\ref{weak_dg}) which is approximated in the discrete system (\ref{semidiscrete}) by $\sum_{m=1}^3 S_{k,x_m} F\left(U_h^{\text{cub}}\right)$
the threads $t_0, \dots, t_{N_{\text{cub}}-1}$ evaluate the flux function $F$ on $U$ in parallel. Finally, the multiplication with the three stiffness matrices $S_{x_1}, S_{x_2}, S_{x_3}$ is treated as described above. However, this matrix vector product can not be handled as one thread per output value, since $S_{x_i}$ has more columns than rows in general. Thus, the threads are distributed over several matrix columns, which are then handled in parallel to maintain the occupancy of the streaming processors high.\\
For this implementation, we were inspired by the MIDG code \cite{midg}.
MIDG is a very lightweight discontinuous Galerkin solver for Maxwell's equations designed to run on multiple GPUs.
A mesh partitioning with ParMetis is involved which minimizes the cuts between processes and thus reduces communication.
The time discretization is realized by a low storage Runge-Kutta scheme as described in section \ref{dg_method}.
For the spatial discretization, the nodal DG scheme as described in \cite{nodalDGbook} is applied.
We mostly adopted the parallel framework of this code including the mesh partitioning and the MPI communications.
However, the quadrature free approach therein does not support curved elements and leads to errors when dealing with nonlinear PDEs like the Euler equations (cf.\ section \ref{results}).
We thus had to apply a computationally more expensive, quadrature based scheme as introduced in section \ref{dg_method}.
\begin{table}[h!]
\begin{center}
\begin{tabular}{r | r | r | r | r | r}
$p$ & CPU & GPU & GPU (float) & speedup & speedup (float)\\ \hline
4 & $278.92$ s & $54.16$ s & $15.27$ s & $5.15$ & $18.27$\\
3 & $92.88$ s & $21.34$ s & $7.96$ s & $4.35$ & $11.67$ \\
2 & $34.01$ s & $8.40$ s & $3.52$ s & $4.05$ & $9.66$ \\
\end{tabular}
\end{center}
\caption{Run time of the time stepping loop, 200 iterations on 25956 cells}
\label{time_table}
\end{table}\\
An overview of our algorithm is given in figure \ref{flowchart}. We work on a compute cluster - shared or distributed memory - on which each node is equipped with one GPU.
Here, the executions on compute node $j$ are described and it can be seen that after the interpolation to the quadrature nodes on the element faces the communication with other processes takes place.
For that purpose, the values on processor cuts are downloaded from the GPU and distributed to other compute nodes. In order to save time, this is done asynchronously, while the volume integration is executed.\\
In our test setting, we have a shared memory system equipped with $8$ Intel Xeon E5620 CPU cores and $8$ Nvidia Tesla M2050 GPUs.
For measurements of the GPU accelerated code against a conventional CPU code, we run a test on one GPU and compare it to a CPU code where the matrix vector multiplication is implemented using BLAS routines. Since the underlying CPU has four cores, we use an OPENMP parallelization to achieve full utilization of its capacity. For that purpose we hint the compiler to apply a parallelization of the k-loop over all elements (cf.\ equation (\ref{weak_dg})). In the test setting of table \ref{time_table} we achieved a run time on one CPU core of $1016.36$ s for $p=4$ elements. Compared to $278.92$ s for the parallel version this reflects the theoretical speedup on this hardware architecture. Furthermore, table \ref{time_table} shows the run times for $200$ Runge-Kutta iterations on a mesh with $25956$ cells for the CPU and GPU code. Here, the most interesting polynomial degrees with respect to CFD are inspected. For the GPU run times we distinguish a single precision (float) and a double precision implementation. The last two columns of table \ref{time_table} show the gained speedup compared to the CPU code. We observe that the speedup increases with the polynomial order. This stems from the higher number of quadrature points in one element leading to a better utilization of the available CUDA threads. Overall, on the most interesting level for this paper of $p=4$ elements, speedup factors of $5$ up to $18$ for single and double precision are achieved.

\end{section}

\begin{section}{Numerical Results}\label{results}

As a special case of equation (\ref{conservation_law}) we consider the Euler equations of gas-dynamics in three spatial dimensions.
These equations describe the motion of an inviscid fluid without heat conduction
\begin{equation*}
	\begin{bmatrix}
		\rho \\
		\rho u \\
		\rho v \\
		\rho w \\
		\rho E \\
	\end{bmatrix}_t
	+
	\begin{bmatrix}
		\rho u \\
		\rho u^2 + p \\
		\rho u v \\
		\rho u w \\
		u ( \rho E + p )\\
	\end{bmatrix}_x
	+
	\begin{bmatrix}
		\rho v \\
		\rho u v \\
		\rho v^2 + p \\
		\rho v w \\
		v ( \rho E + p )\\
	\end{bmatrix}_y
	+
	\begin{bmatrix}
		\rho w \\
		\rho u w \\
		\rho v w \\
		\rho w^2 + p \\
		w ( \rho E + p )\\
	\end{bmatrix}_z
	=
	\begin{bmatrix}
		0 \\
		0 \\
		0 \\
		0 \\
		0 \\
	\end{bmatrix}.
\end{equation*}
Here, the unknown function reads as
\begin{equation*}
U: \mathbb{R} \times \mathbb{R}^3 \to \mathbb{R}^5 \, , \quad U(t,\boldsymbol{x}) = \left(\rho(t,\boldsymbol{x}),\, \rho u(t,\boldsymbol{x}),\, \rho v(t,\boldsymbol{x}),\, \rho w(t,\boldsymbol{x}),\, \rho E(t,\boldsymbol{x})\right)^T
\end{equation*}
where  $\rho$ denotes the density, $u, v, w$ the fluid velocities in the three space directions and $E$ the total energy.
Finally, the system is completed by the perfect gas law for the pressure
\begin{equation*}
	p = (\gamma - 1) \left( \rho E - \frac{u^2 + v^2 + w^2}{2\rho} \right)
\end{equation*}
where $\gamma$ is the adiabatic index of the fluid \cite{blazek}.\\
%It should be mentioned that the flux vectors are rational functions in terms of the conservative variables, which forces us to apply very high order quadrature rules.\\
Our first test case is a subsonic flow past a sphere at $M_{\infty} = 0.38$ which was examined together with the discontinuous Galerkin method in \cite{bassi_rebay_euler}.
This test case is very attractive, since the surface curvature is straightforward. Nevertheless, we apply our curvature approach to obtain a smooth volume mesh.
Let $r$ denote the radius and $\boldsymbol{x_0}$ the centroid of the sphere. Then the boundary conditions in equation (\ref{boundary}) for the mesh deformation  are given by
\begin{equation*}
g(\boldsymbol{x}) = \boldsymbol{x_0} + \frac{r}{\left\Vert \boldsymbol{x} - \boldsymbol{x_0} \right\Vert }\boldsymbol{x} - \boldsymbol{x}.
\end{equation*}
\begin{figure}
\subfigure[Disturbed Solution on straight sided elements]{\includegraphics[width=0.49\textwidth]{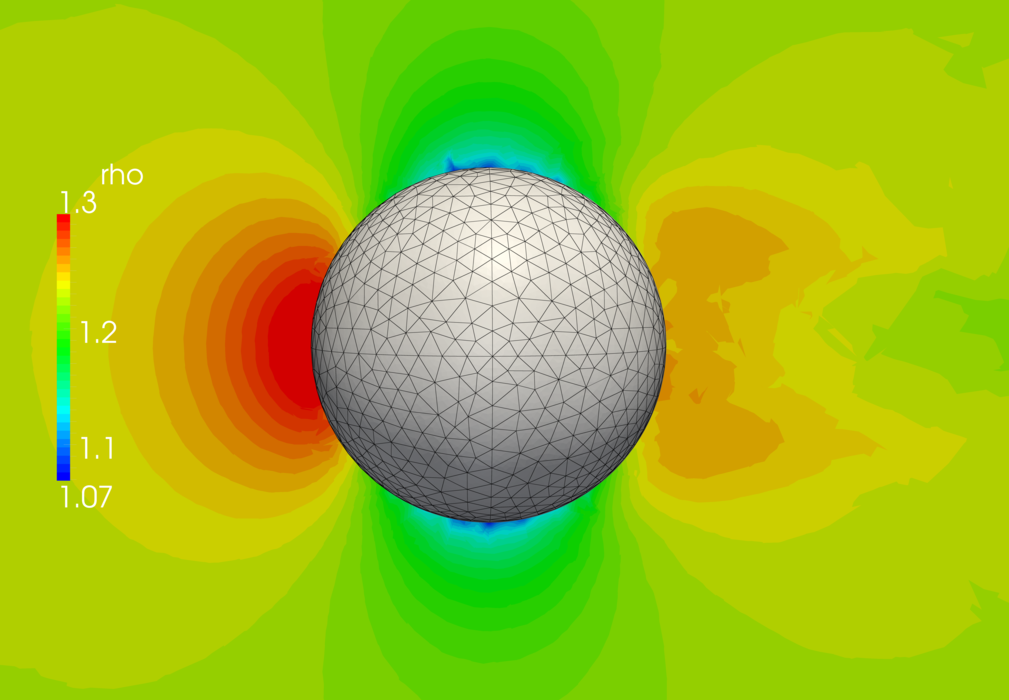}}\hfill
\subfigure[Isoparametric curved elements]{\includegraphics[width=0.49\textwidth]{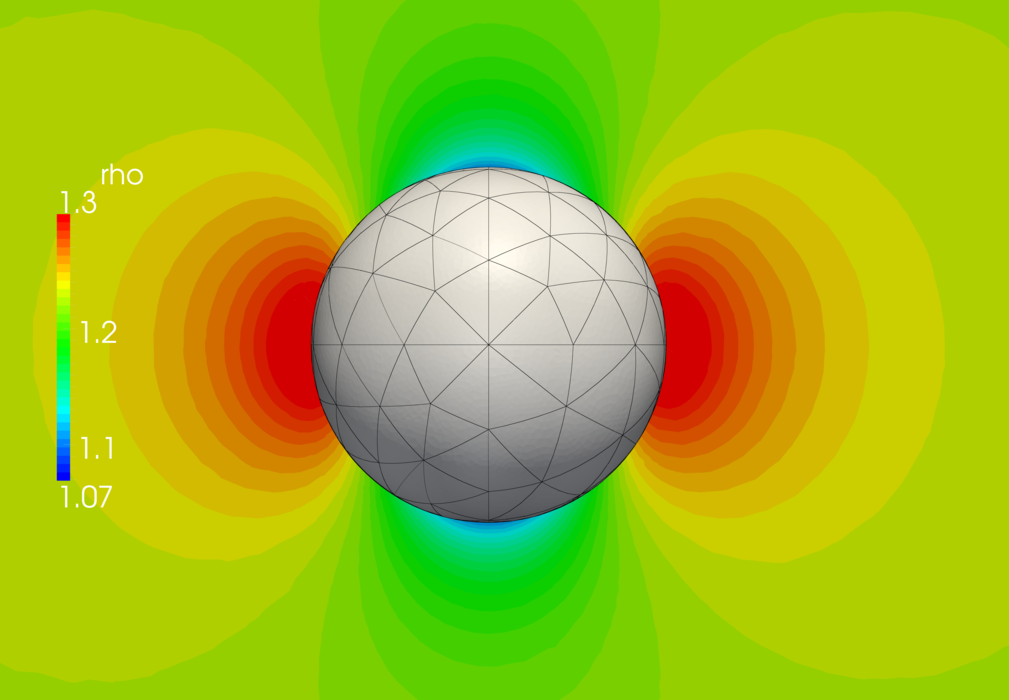}}
\caption{Density distribution of a subsonic flow past a sphere with $p=4$ basis functions}
\label{sphere_results}
\end{figure}\\
As shown in figure \ref{sphere}, we start with a straight sided grid and solve the linear elasticity equation on the surrounding volume grid with $p=4$ basis functions.
Thus, we obtain a smoothly curved volume discretization, which we then use as input for the DG Euler solver also based on $p=4$ basis functions.
We note that best results are obtained when the order of the linear elasticity solution matches the order of the DG solver. In this case the deformation lies within the range of the DG basis functions and can be resolved properly.\\ 
Figure \ref{sphere_results} shows the effect of curved boundaries on the solution quality. On the left hand side the solution is disturbed due to straight sided elements in the discretization mesh. This leads to small yet problematic kinks between the elements on the surface and results in a non-physical solution. In contrast, the computations on the right hand side were performed on a curved grid, which is in addition coarser. In this case the linear elasticity deformation was calculated using polynomials of order four. 
\\The next test case is the NACA0012 symmetric airfoil.
Although the geometry is only two dimensional, we stretch the profile into the third coordinate direction in order to apply the same solver as for the other test cases.
Here, the effect of the boundary deformation on the volume mesh can be visualized as figure \ref{naca} shows.\\
In the NACA0012 situation we consider two flow conditions. First, a pure subsonic flow with no angle of attack and Mach number $M_{\infty}=0.4$ (figure \ref{naca_results}a). 
Figure \ref{naca_results}b shows a $M_{\infty}=0.8$ transonic flow with angle of attack $\alpha = 1.25^\circ$. In both pictures, the density distribution is plotted.
Since we are dealing with discontinuities in this situation, we have to add artificial viscosity to the system of equations in order to ensure stability of the method.
For that goal, we proceed as described in section \ref{dg_method} and apply a shock detector to select the troubled cells.
\\In the following we consider the Onera M6 wing. This test case is well known and examined with a variety of numerical methods.
The original experiment was defined in \cite{onera} at Mach number $M_{\infty} = 0.8395$ and angle of attack $\alpha  = 3.06^{\circ}$.
Again, we start with a straight-sided, tetrahedral mesh and a NURB representation of the ONERA M6 geometry and use the deformation approach to obtain the curved mesh (figure \ref{onera}).
Figure \ref{onera_results} shows the density distribution on the upper surface of the airfoil.
\begin{figure}
\subfigure[Subsonic flow, $M_{\infty}=0.4, \alpha = 0^\circ$]{\includegraphics[width=0.49\textwidth]{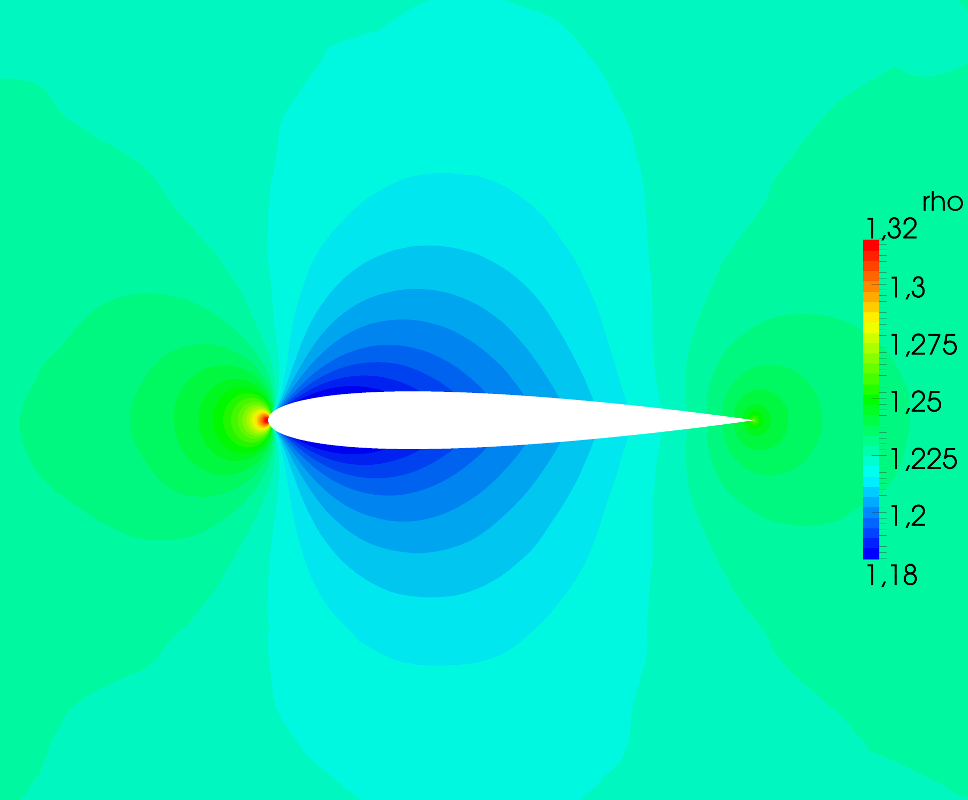}}\hfill
\subfigure[Transsonic flow, $M_{\infty}=0.8, \alpha = 1.25^\circ$]{\includegraphics[width=0.49\textwidth]{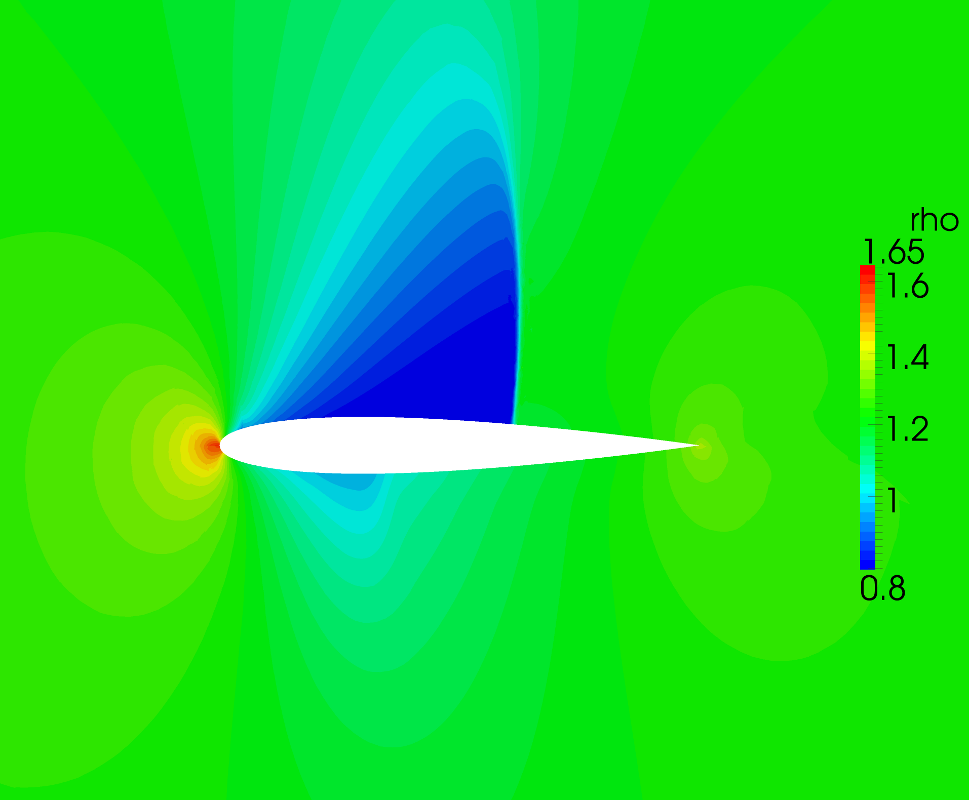}}
\caption{Flow past NACA0012 airfoil with $p=4$ elements}
\label{naca_results}
\end{figure}
Also in this case, artificial viscosity is applied to the system to deal with the shocks on the upper surface of the airfoil.
In both the NACA0012 and ONERA M6 situation we set the parameter for the viscosity $\kappa = 4$ and $\epsilon_0 = 0.3$.\\
For this computations we utilize a so-called $p$-refinement as a acceleration  technique. This means that we start the time stepping with low order polynomials of degree $p_1$, e.g.\ $p_1 = 1$, iterate until the Runge-Kutta update is below a specified tolerance and then switch to a richer polynomial space of degree $p_2 = p_1 +1$. This process is repeated until the desired order is reached. The transfer operation between the coarser and the finer polynomial space is straightforward, since we are dealing with a hierarchical set of basis functions as introduced in section \ref{dg_method}. Thus, the solution on level $p_i$ can be directly embedded into the $p_{i+1}$ space, where the higher order polynomials are initially weighted with zeros.\\
This procedure offers two advantages in terms of acceleration. Firstly, computations are much cheaper on a coarser level of $p$, since the number of unknown values is smaller. Furthermore, the cubature and quadrature formulas do not have to be as precise as for higher order polynomials, which leads to smaller local operators. Secondly, the timestep for the Runge-Kutta time integration scales with $\mathcal{O}(p^{-2})$, c.f. \cite{nodalDGbook}. This enables a much larger timestep for smaller values of $p$, which reduces the number of iterations needed to reach the steady state solution.\\
For the computations shown in this work, we started with $p=2$ basis functions and subsequently refined to $p=3$ and $p=4$. We found out that this is a convenient choice, since we are dealing with relatively coarse grids, where computations with $p=1$ or even $p=0$ do not lead to suitable results.\\ 
\begin{table}[h!]
\begin{minipage}{0.5\textwidth}
\begin{tabular}{r | r | r | r}
$p$ & iterations & timestep & time consumed\\ \hline
2 & $8000$ & 6e-5& $105.8$ s\\
3 & $10000$& 2e-5& $196.5$ s\\
4 & $20000$& 5e-6& $650.0$ s\\ \hline
 & $38000$& & $952.3$ s
\end{tabular}
\end{minipage}
\begin{minipage}{0.5\textwidth}
\begin{tabular}{r | r | r | r}
$p$ & iterations & timestep & time consumed\\ \hline
2 &  & & \\
3 & & & \\
4 & $61000$& 5e-6& $1982$ s\\ \hline
& $61000$& & $1982$ s
\end{tabular}
\end{minipage}
\caption{Algorithm execution for the sphere - with and without p-refinement}
\label{table_p_refine}
\end{table}
Table \ref{table_p_refine} shows two different algorithm executions for the sphere test case. The underlying discretization mesh consists of 25956 tetrahedra and 4457 of those are curved elements. Moreover, on the finest level, each element contains 35 collocation nodes, 70 cubature nodes and $4\cdot 16$ surface quadrature nodes. On the left hand side, the simulation was started with $p=2$ basis functions and then subsequently refined to $p=3$ and $p=4$ functions. Here, the number of iterations on  level $2$ and $3$ was chosen empirically. In contrast, on the right hand side, the simulation was run using $p=4$ basis functions only. Both executions terminated when the infinity norm of the residual for $p=4$ was below a tolerance of 1e-9. For this particular case, the execution time is more then halved by the p-refinement procedure. The number of iterations on each level may look surprising, since it is multiple of $1000$. This stems from the fact that the calculation of a norm is an expensive operation in parallel. In contrast, a few extra Runge-Kutta evaluations are relatively cheap compared to the parallel reduction for the norm. Thus, we evaluate the norm of the residual only after multiples of $1000$ iterations.
\begin{figure}
\includegraphics[width=1.0\textwidth]{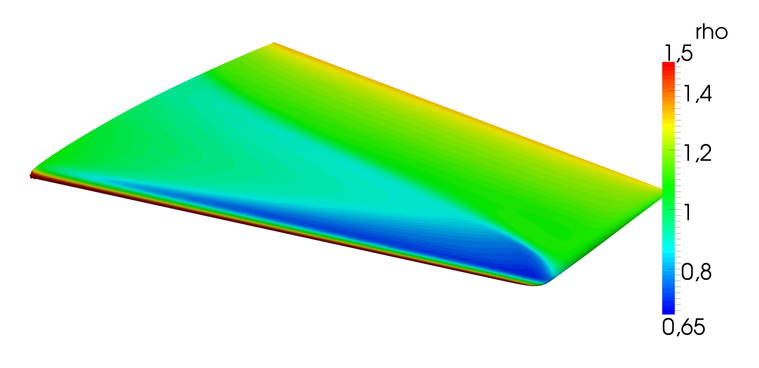}
\caption{Flow over Onera M6 airfoil with $p=4$ elements}
\label{onera_results}
\end{figure}
\end{section}

\begin{section}{Conclusion}\label{conclusion}
In this work, we have presented a curved mesh generation approach.
This approach is based on the linear elasticity deformation of an initial grid with straight sided elements until it meets the desired curved boundary.
We demonstrated that the boundary approximation can be of arbitrary high order reducing the error induced into the discontinuous Galerkin discretization due to under-resolved geometries.
The solution of the linear elasticity equations with a high order finite element method seems computationally very expensive at a first glance. However, this has to be done only once.
We precompute the curvature and the DG solver only needs to evaluate the FEM-solution at the desired nodes.\\
In a second step, we showed how these curved meshes are embedded into a GPU based parallel DG solution of the Euler equations of gas-dynamics.
Furthermore, the performance of this massively parallel GPU code was tested where we gained a speedup of up to $18$ compared to the serial version of this code.
We have chosen some challenging test cases including transonic flows leading to discontinuities which were treated with an artificial viscosity approach.\\
Although the combination of an explicit Runge-Kutta time integration with diffusion terms leads to very small time steps we have seen that the outstanding parallel performance of the GPU overcomes this issue.
And we believe that this fact justifies the additional implementation effort for the GPU code.
\end{section}

\section*{Acknowledgements}
Our research is supported by BMBF (Bundesministerium für Bildung und Forschung) within the collaborative project DGHPOPT.

\bibliographystyle{unsrt}
\bibliography{dg_paper}

\end{document}